\documentclass[12pt]{article}
\usepackage{amsmath,amssymb,amsfonts,graphicx,wrapfig}

\newtheorem{theorem}{Theorem}
\newtheorem{lemma}{Lemma}

\begin{document}

\begin{center}

{\LARGE
The distribution of the maximum number of common neighbors in the random graph
}

{\Large

I.V. Rodionov\footnote{Institute for Information Transmission Problems (Kharkevich Institute) of the Russian Academy of Sciences, Moscow, Russian Federation.\\vecsell@gmail.com}, M.E. Zhukovskii\footnote{Moscow Institute of Physics and Technology (State University), laboratory of advanced combinatorics and network applications, Dolgoprodny, Moscow Region, Russian Federation; Adyghe State University, Caucasus mathematical center, Maykop, Republic of Adygea, Russian Federation; The Russian Presidential Academy of National Economy and Public Administration, Moscow, Russian Federation.\\zhukmax@gmail.com}

}

\vspace{0.5cm}

Abstract\\

\end{center}

Let $\Delta_{k;n}$ be the maximum number of common neighbors of a set of $k$ vertices in $G(n,p)$. In this paper, we find $a_n$ and $\sigma_n$ such that $\frac{\Delta_{k;n}-a_n}{\sigma_n}$ converges in distribution to a random variable having the standard Gumbel distribution.

\vspace{0.3cm}

{\bf Keywords:} binomial random graph, maximum degree, common neighbors, Gumbel distribution

\vspace{0.5cm}

\section{Introduction}
\label{Intro}

In 1980~\cite{Bol_Degree}, B. Bollob\'{a}s studied the asymptotical behavior of the maximum degree $\Delta_n$ of the binomial random graph $G(n,p)$ (\cite{Bollobas,Janson}) for fixed $p\in(0,1)$. The main result of that paper is the following. Let, for $n\in\mathbb{N}$, $n\geq 2$,
\begin{equation}
a_n=pn+\sqrt{2p(1-p)n\ln n}\left(1-\frac{\ln\ln n}{4\ln n}-\frac{\ln(2\sqrt{\pi})}{2\ln n}\right),
\label{a_n}
\end{equation}
\begin{equation}
\sigma_n=\sqrt{\frac{p(1-p)n}{2\ln n}}.
\label{sigma_n}
\end{equation}
Then the shifted random variable $\frac{\Delta_n-a_n}{\sigma_n}$ converges in distribution to a random variable having the standard Gumbel distribution. Further in his paper, Bollob\'{a}s considered the $m$-th highest degree $\Delta_n^m$ of $G(n,p)$ (in particular, $\Delta^1_n=\Delta_n)$ and proved a similar result: for every $y\in\mathbb{R}$,
$$
{\sf P}\left(\frac{\Delta_n^m-a_n}{\sigma_n}\leq y\right)\to e^{-e^{-y}}\sum_{j=0}^{m-1}\frac{e^{-jy}}{j!}\text{ as }n\to\infty.
$$
Note that elements $\xi_i$ of a sequence of degrees of $G(n,p)$ have the binomial distribution with the parameters $n-1,p$. For sequences of series of independent binomial random variables, the asymptotical distribution of maximums was studied by S. Nadarajah and K. Mitov in 2002~\cite{Independent}. Their result states that the maximum $D_n$ of $n$ independent binomial random variables $\xi_1^{(n)},\ldots,\xi_n^{(n)},$ where the upper index denotes the number of series, with parameters $M = M(n)=\omega(\ln^3 n)$ and $p=\mathrm{const}$ (here and below we use the usual notation $g_1= \omega(g_2)$ for two sequences $g_1(n)$ and $g_2(n)$ such that $g_2 = o(g_1)$) obeys the following asymptotical law: for $y\in\mathbb{R}$,
\begin{equation}
{\sf P}\left(D_n\leq pM+\sqrt{2p(1-p)M\ln n}\left[1-\frac{\ln\ln n}{4\ln n}-\frac{\ln(2\sqrt{\pi})}{2\ln n}+
\frac{y}{2\ln n}\right]\right)\to e^{-e^{-y}}\text{ as }n\to\infty.
\label{approaches_Gumbel}
\end{equation}
It is easy to see that, for $M=n-1$, the result of Nadarajah and Mitov gives the same normalization functions $a_n$ and $\sigma_n$ and the same asymptotical distribution as the result of Bollob\'{a}s (but for {\it dependent} random variables). Clearly, a simple substitution gives slightly different functions $\tilde a_n$ and $\tilde\sigma_n$. However, the convergence of ${\sf P}(\Delta_n\leq\tilde a_n+y\tilde\sigma_n)$ implies the convergence of ${\sf P}(\Delta_n\leq a_n+y\sigma_n)$ to the same limit since $a_n=\tilde a_n+O(1)$, $\sigma_n=\tilde\sigma_n+O((n\ln n)^{-1})$, and $\sigma_n=\Theta(\sqrt{n}/\sqrt{\ln n})$ (the relation $g_1= \Theta(g_2)$ for two sequences $g_1(n)$ and $g_2(n)$ means that $g_1 = O(g_2)$ and $g_2 = O(g_1)$).\\

Results of such kind belong to the extreme value theory. A general result of this theory, the Fisher--Tippett--Gnedenko theorem~\cite{Gnedenko} (see also~\cite{David}, page 205) states the following. Let $\xi_1,\xi_2,\ldots$ be independent and identically-distributed random variables, $\xi^{(n)}=\max\{\xi_1,\ldots,\xi_n\}$. If there exist $a_n\in\mathbb{R}$, $\sigma_n>0$ and a non degenerate distribution $F$ such that $\frac{\xi^{(n)}-a_n}{\sigma_n}$ converges to an $\eta\sim F$ in distribution, then $F$ belongs to either the Gumbel, the Fr\'{e}chet or the Weibull family. In this way, three extremal types of distributions of $\xi_i$ are distinguished. A comprehensive account of necessary and sufficient conditions for a distribution to belong to one of the extremal types is given in~\cite{Leadbetter}. For further results of the extreme value theory and its applications see, e.g., \cite{Beirlant,Haan}. However, mentioned results can not be applied to the above case (when $\xi_1,\ldots,\xi_n$ are identically distributed but their distribution depends on $n$). Such {\it triangular arrays} of random variables were studied, e.g., in~\cite{Anderson,Sielenou,Independent}.

The result of Nadarajah and Mitov easily follows from a theorem about large deviations for Binomial random variables and certain properties of the normal distribution function. But how can~(\ref{approaches_Gumbel}) be obtained for dependent random variables $\xi_1,\ldots,\xi_n$? Bollob\'{a}s, in his proof, introduced the random variable $X$ being the number of vertices having degree greater than $y\sigma_n+a_n$, where $a_n$ and $\sigma_n$ are defined in (\ref{a_n}) and (\ref{sigma_n}) respectively. The result follows from the fact that, for every $j\in\mathbb{N}$, the $j$-th factorial moment of $X$ converges to the $j$-th factorial moment of a Poisson random variable with the parameter $e^{-y}$. Obviously, the same idea may be used to prove the result of Nadarajah and Mitov and the respective result for identically distributed (not depending on $n$) independent random variables under certain conditions (see~\cite{Leadbetter}, Chapter 2).\\

Let $k$ be an arbitrary positive integer. In our paper, we solve the problem of finding an asymptotical distribution (precisely, the normalizing sequences $a_{k;n}$, $\sigma_{k;n}$ and an extremal type) of the maximum $\Delta_{k;n}$ of  the number of common neighbors of $k$ vertices in $G(n,p)$ (the case $k=1$ is already solved by Bollob\'{a}s since $\Delta_{1;n}=\Delta_n$). Notice that, in our paper, we consider not only constant $p$ but $p$ depending on $n$. In~\cite{Ivchenko}, it is proven that the result of Bollob\'{a}s holds true for $p=p(n)\to 0$ as $n\to\infty$ such that $\frac{pn}{\ln^3 n}\to\infty$ as $n\to\infty$ (convergence to the Gumbel distribution but for {\it other} functions $a_n,\sigma_n$ holds true even if $\frac{pn}{\ln n}\to\infty$, and the latter condition is optimal). In our paper, we obtain a similar result for $\Delta_{k;n}$ in the same most general settings.

For $k\geq 2$, this problem differs a lot from both mentioned problems (the case of independent binomial random variables and degrees of the random graph). The main difference is that the variance of an analogue of the random variable $X$ defined above may approach infinity (e.g., this happens when $k=2$ and $p>1/2$). In particular, this fact makes it impossible to apply the method of factorial moments directly.

More formally, let $v_1,\ldots,v_k\in[n]:=\{1,\ldots,n\}$ be distinct vertices of $G(n,p)$. Let $N_n(v_1,\ldots,v_k) \subseteq [n]\backslash\{v_1,\ldots,v_k\}$ be the set of all common neighbors of $v_1,\ldots,v_k$ in $G(n,p)$ ($u\in N_n(v_1,\ldots,v_k)$ if and only if, for every $i\in[k]$, $u\sim v_i$, i.e., $u$ is adjacent to $v_i$).  Set
$$
\Delta_{k;n}=\max_{v_1,\ldots,v_k}|N_n(v_1,\ldots,v_k)|,
$$
where the maximum is over all distinct vertices $v_1,\ldots,v_k\in[n]$. Moreover, let $\Delta^m_{k;n}$ be the $m$-th highest value among $|N_n(v_1,\ldots,v_k)|$ (in particular, $\Delta_{k;n}=\Delta^1_{k;n}$). The main result of our paper is given below.

\begin{theorem}
Fix $y\in\mathbb{R}$ and $k,m\in \mathbb{N}.$ Assume $p = p(n)\in(0,1)$ is such that
\begin{equation}
p^k\gg\frac{\ln^3 n}{n},\quad
1-p\gg\sqrt{\frac{\ln\ln n}{\ln n}}.
\label{cond}
\end{equation}
Let
$$
a_{k;n}=np^k+\sqrt{2kp^k(1-p^k)n\ln n}\left(1-\frac{\ln[k!]}{2k\ln n} - \frac{\ln[4\pi k \ln n]}{4k\ln n}\right),
$$
$$
\sigma_{k;n}=\sqrt{\frac{p^k(1-p^k)n}{2k\ln n}}.
$$
Then
$$
 {\sf P}\left(\frac{\Delta^m_{k;n}-a_{k;n}}{\sigma_{k;n}}\leq y\right)\to e^{-e^{-y}}\sum_{i=0}^{m-1}\frac{e^{-yi}}{i!}\quad\text{as }n\to\infty.
$$
\label{main}
\end{theorem}

{\it Remark}. The second condition may be strengthened: $1-p\geq\sqrt{\frac{32\ln\ln n}{k^3\ln n}}(1+o(1))$. For such $p$, the same techniques work. However, we give a proof in a weaker form to avoid some annoying technical details.\\

Note that $\Delta_{k;n}$ is the maximum over ${n\choose k}$ binomial random variables with parameters $n-k,p^k$. Therefore, our result duplicates the statement~(\ref{approaches_Gumbel}) (but for this special case of dependent random variables). This motivates the following question. How strong can be dependencies between binomial random variables until~(\ref{approaches_Gumbel}) fails? A partial answer on this question and other further questions are given in Section~\ref{Further}.\\


While Theorem~\ref{main} is a natural extension of the result of Bollob\'{a}s, it is also motivated by the study of {\it extension counts}. This study was initiated by Spencer in~\cite{Spencer_ext_counting}. He proved that, given a {\it strictly balanced grounded} pair of graphs $(G,H)$ (where $G$ has $k$ vertices) and $\varepsilon>0$, there exists $C=C(\varepsilon)>0$ such that if $\mu>C\ln n$, then
\begin{equation}
{\sf P}\left(\max_{v_1,\ldots,v_k} |X(v_1,\ldots,v_k)-\mu|<\varepsilon\mu\right)\to 1\quad\text{as }n\to\infty.
\label{ext_law}
\end{equation}
Here, $X(v_1,\ldots,v_k)$ is the number of {\it $(G,H)$-extensions} of the tuple $(v_1,\ldots,v_k)$ and $\mu={\sf E}X(1,\ldots,k)$.

This result was recently strengthened by \v{S}ileikis and Warnke in~\cite{Warnke_extensions}. They proved that there exist constants $c,C,\alpha>0$ such that, for all $p=p(n)\in[0,1]$ and $\varepsilon=\varepsilon(n)\in[n^{-\alpha},1]$, the limit probability in~(\ref{ext_law}) equals 0 if $\mu<\frac{c}{\varepsilon^2}\ln n$, and equals 1 if $\mu>\frac{C}{\varepsilon^2}\ln n$. Theorem~\ref{main} implies, in particular, that, for $p$ such that $p\gg\sqrt{\frac{\ln\ln n}{\ln n}}$ and $1-p\gg\sqrt{\frac{\ln\ln n}{\ln n}}$ (the first restriction on $p$ appears because we also need an asymptotical distribution of the minimum number of $(G,H)$-extensions, and it follows from Theorem~\ref{main} by considering the complement of $G(n,p)$ which is distributed as $G(n,1-p)$)  and $H$ having one more vertex than $G$ adjacent to all the vertices of $G$, the threshold for the concentration result~(\ref{ext_law}) is fully determined: the result of \v{S}ileikis and Warnke is true for any $c<2k(1-p^k)\left(1-(1-\delta)\frac{\ln\ln n}{2k\ln n}\right)$, $C>2k(1-p^k)\left(1-(1+\delta)\frac{\ln\ln n}{2k\ln n}\right)$ and any $\alpha>0$. Notice that $c$ is bounded away from 0 only when $p$ is bounded away from 1. However, it does not contradict the result of \v{S}ileikis and Warnke since our conclusion becomes non-trivial only when $\varepsilon$ is close to $n^{-1/2}$ since, in our settings, $\mu=n^{1-o(1)}$. We hope that our methods can be further applied to get similar results for other pairs of $G$ and $H$.\\


The rest of the paper is organized as follows. In Section~\ref{Proof}, we give a proof of Theorem~\ref{main}. Section~\ref{Further} is devoted to some discussions of our method and its possible applications to more general questions.

\section{Proof of Theorem~\ref{main}}
\label{Proof}

Consider the random variables $X=X_n^k$ being the number of $k$-sets of vertices having more than
$$
b=b_{k;n}(y):=a_{k;n}+y\sigma_{k;n}
$$
common neighbors. The reason why $\mathrm{Var}X$ may approach infinity is that the major contribution to the variance is made by those $k$-sets that have proper subsets with large number of common neighbors. A.a.s., in $G(n,p)$ there are no such $k$-sets (for details, see Section~\ref{Proof_Moments}). So, we are able to exclude such $k$-sets from $X$ (in what follows, we denote the shifted random variable by $\tilde X$).  In Section~\ref{Proof_Moments}, we estimate the expectation and the second moment of $\tilde X$.

The second problem we face is that a direct implementation of the approach of Bollob\'{a}s (i.e., estimation of all the factorial moments of $\tilde X$) requires heavy computations. Fortunately, we have shown that it is enough to prove that ${\sf E}\tilde X(\tilde X-1)\sim({\sf E}X)^2$. This observation follows from the Janson-type inequality that we prove in Section~\ref{Proof_Janson} (in Section~\ref{Further}, we state it in a more general form and describe its possible applications to other problems related to extreme value theory). It is worth mentioning that the result of Bollob\'{a}s would follow directly from the Janson inequality for general upsets~(\cite{Riordan}, Theorem 1, Inequality (3)), if, for distinct vertices $u,v$, the events $A_u,A_v$ of having more than $b_{1;n}(y)$ neighbors of $u$ and having more than $b_{1;n}(y)$ neighbors of $v$ respectively were independent (indeed, for a fixed vertex $v$, the property $A_v$ is {\it increasing}, i.e. is the upset). If the premise held, estimating the first moment would be enough to get the asymptotics of ${\sf P}(X_n^1=0)$. Unfortunately, any two $A_u,A_v$ are dependent, and so, the mentioned Janson inequality implies ${\sf P}(X_n^1=0)\leq e^{-e^{-y}+\frac{1}{2}e^{-2y}+o(1)}$.

Nevertheless, we prove Janson-type bounds for ${\sf P}(\tilde X=0)$ for all $k$ (Section~\ref{Proof_Janson}, Lemma~\ref{Janson-type}). It is easy to see that, for $k=1$, the same result holds true for the original (not shifted) random variable $X_n^1$: $e^{-{\sf E}X_n^1+o(1)}\leq{\sf P}(X_n^1=0)\leq e^{-{\sf E}X_n^1+o\left(e^{{\sf E}X_n^1}\right)}$, and it immediately implies the result of Bollob\'{a}s on $\Delta_n$.\\ 



\subsection{Notations and approximations of binomial distributions}

In our proofs, we will frequently use the classical De Moivre--Laplace central limit theorem. For the sake of convenience, below, we verify that it gives (in our settings~(\ref{cond})) an approximation error which is sufficient for our goals.

Given $q=q(n)\in(0,1)$ and $x=x(n)\in\mathbb{R}$, assume that
\begin{equation} \label{assumptions_moivre}
\frac{\ln^3 n}{nq(1-q)}\to 0,\quad x^3\sqrt{\frac{\ln^3 n}{nq(1-q)}}\to 0\quad\text{ as }n\to\infty.
\end{equation}
Below, we find asymptotics for ${\sf P}(\xi_{n,q} = a)$ and ${\sf P}(\xi_{n,q} > a)$ where
$$
a = \left\lfloor nq + x\sqrt{q(1-q) n \ln n}\right\rfloor
$$
and $\xi_{n,q}$ is a binomial random variable with parameters $n,q$. 

Using Stirling's formula, the relation $\ln(1+z)=z-\frac{z^2}{2}+O(z^3)$ and the definition of $a,$ we get
$$
{\sf P}(\xi_{n,q} = a) = \binom{n}{a} q^a (1 - q)^{n-a} = \frac{1 + O(1/a)+O(1/(n-a))}{\sqrt{2\pi n}\sqrt{a/n} \sqrt{1 - a/n}}\frac{q^a (1-q)^{n-a}}{(a/n)^a(1 - a/n)^{n-a}}
$$
$$ = \frac{1 + O\left(x\sqrt{\ln n/(nq(1-q))}+1/a+1/(n-a)\right)}{\sqrt{2\pi n q (1 - q)}} \exp\left\{-n \ln \left[\left(\frac{a/n}{q}\right)^{a/n} \left(\frac{1-a/n}{1-q}\right)^{1-a/n}\right] \right\}
$$
$$ \sim \frac{1}{\sqrt{2\pi n q (1 - q)}}\exp\left\{-nq \left(1 + x\sqrt{\frac{(1-q)\ln n}{q n}}\right) \ln\left(1 + x\sqrt{\frac{(1-q)\ln n}{q n}}\right) \right.
$$
$$
\left. - n(1-q)\left(1 - x\sqrt{\frac{q\ln n}{(1-q) n}}\right)\ln \left(1 - x\sqrt{\frac{q\ln n}{(1-q) n}}\right)\right\}
$$
\begin{equation}
= \frac{1}{\sqrt{2\pi n q (1 - q)}} \exp\left\{-\frac{x^2 \ln n}{2} + O\left(x^3\sqrt{\frac{\ln^3 n}{nq(1-q)}}\right)\right\}\sim\frac{1}{\sqrt{2\pi nq(1-q)}}\exp\left\{-\frac{x^2 \ln n}{2}\right\}.
\label{demoivre}
\end{equation}

Assume that (in addition to the above conditions)
\begin{equation} \label{assumptions_moivre1}
x\geq 0,\quad x\sqrt{\ln n}\to\infty,\quad \frac{x^3}{(1-q)^2}\sqrt{\frac{\ln^3 n}{n}}\to 0\quad\text{ as }n\to\infty.
\end{equation}
Denote $\hat x=\frac{2x}{\sqrt{1-q}}$. 
Let
$$
\hat a=nq+\hat x\sqrt{q(1-q) n\ln n}, \quad
a^+=a, \, a^-=nq - x\sqrt{q(1-q) n \ln n}.
$$
Then, approximating the sum by an integral, we have
\begin{multline*}
{\sf P}(\xi_{n,q}>a^+) = \\
\sum_{j=a+1}^{\hat a} \left(1 + O\left(\frac{x^3}{(1-q)^2}\sqrt{\frac{\ln^3 n}{nq}}\right)\right)
\frac{\exp\left(-(j-nq)^2/[2nq(1-q)]\right)}{\sqrt{2\pi n q(1-q)}}+{\sf P}(\xi_{n,q}>\hat a)=\\
 \int_{x\sqrt{\ln n}}^{\hat x\sqrt{\ln n}}\frac{1}{\sqrt{2\pi}}e^{-\frac{t^2}{2}}dt \left(1 + O\left(\frac{x^3}{(1-q)^2}\sqrt{\frac{\ln^3 n}{nq}}\right)\right)+{\sf P}(\xi_{n,q}>\hat a) =\\
 \left(1 - \Phi\left(x\sqrt{\ln n}\right)\right)(1 + o(1))-\left(1 - \Phi\left(\hat x\sqrt{\ln n}\right)\right)(1+o(1))+{\sf P}(\xi_{n,q}>\hat a),
\end{multline*}
where $\Phi(t)=\int_{-\infty}^t \frac{1}{\sqrt{2\pi}}e^{-z^2/2}dz$. It is also worth to be noted that the $O$-term inside the integral on the second line of the latter relation converges to $0$ by \eqref{assumptions_moivre} if $q$ is bounded away from 1 and by \eqref{assumptions_moivre1} if $q\to 1.$ Applying the well-known relation
\begin{equation}
1-\Phi(t) \sim \frac{1}{\sqrt{2\pi}t}e^{-t^2/2} \quad\text{ as }t\to\infty
\label{Phi_approx}
\end{equation}
(see ($1^{\prime}$) in \cite{Bol_Degree2}), using $\hat{x} \ge 2x$ and $x\sqrt{\ln n}\to\infty$, we get
\begin{eqnarray*}
{\sf P}(\xi_{n,q}>a^+) =
\frac{\mathrm{exp}\left[-x^2\ln n/2\right]}{x\sqrt{2\pi\ln n}}(1 + o(1))-\frac{\mathrm{exp}\left[-\hat x^2\ln n/2\right]}{\hat x\sqrt{2\pi\ln n}}(1 + o(1))+{\sf P}(\xi_{n,q}>\hat a) =\\
 \frac{\mathrm{exp}\left[-x^2\ln n/2\right]}{x\sqrt{2\pi\ln n}}(1 + o(1)) + {\sf P}(\xi_{n,q}>\hat a).
\end{eqnarray*}
Next, by the Chernoff bound (\cite{Janson}, Theorem 2.1), for large $n$,
$$
 {\sf P}\left(|\xi_{n,q}-nq|>\hat a-nq\right)\leq 2\mathrm{exp}\left[-\frac{\hat x^2 q(1-q) n\ln n}{2\left(nq+\hat x\sqrt{ q(1-q)n\ln n}/3\right)}\right]\leq 2\mathrm{exp}\left[-\frac{4x^2 \ln n}{3}\right].
$$
Since
$$\mathrm{exp}\left[-\frac{4x^2 \ln n}{3}\right] = o\left(\frac{\mathrm{exp}\left[-x^2\ln n/2\right]}{x\sqrt{2\pi\ln n}}\right),$$
we finally get
\begin{equation}
{\sf P}(\xi_{n,q}>a^+)\sim
\frac{\mathrm{exp}\left[-x^2\ln n/2\right]}{x\sqrt{2\pi\ln n}}.
\label{central1}
\end{equation}
In the same way,
$$
{\sf P}(\xi_{n,q}<a^-)\sim\frac{\mathrm{exp}\left[-x^2\ln n/2\right]}{x\sqrt{2\pi\ln n}}.\eqno(10')
$$


In the remaining parts of the proof, we use the following notations. For every $\ell\in[k-1]$, denote
$$
\Gamma_{\ell} = np^{\ell} + \sqrt{2\ell}\sqrt{np^{\ell}(1-p^{\ell})\ln n}.
$$

Let $\tilde X = \tilde X_n^k$ be the number of $k$-sets of vertices $U_k:=\{u_1, \ldots, u_k\}$ such that, for every $\ell\in[k-1]$ and every distinct $i_1,\ldots,i_{\ell}\in[k]$, the following conditions hold:
$$
\left|N_n\left(u_{i_1}, \ldots, u_{i_{\ell}}\right)\right| \leq \Gamma_{\ell} \text{ and } |N_n(u_1,\ldots, u_k)|>b.
$$
Fix a $k$-set $U=\{u_1,\ldots,u_k\}\in{[n]\choose k}$. For $W\subset [n]\setminus U$ and $W_i\subset[n]\setminus\{u_i\}$, $i\in[k]$, consider the following events
$$
B_{U,W}=\{W=N_n(U)\}, \quad B_{U,W_1,\ldots,W_k}=\{W_1=N_n(u_1),\ldots,W_k=N_n(u_k)\}.
$$
Let
$$
B_{U}=\bigvee_{W\subset [n]\setminus U:\,|W|>b}B_{U,W},\quad
\tilde B_U=\bigvee B_{U,W_1,\ldots,W_k},
$$
where the second union is over sets $W_i\subset[n]\setminus\{u_i\}$, $i\in[k]$, such that $|W_1\cap\ldots\cap W_k|>b$ and, for every $\ell\in[k-1]$ and every distinct $i_1,\ldots,i_{\ell}\in[k]$, $\left|W_{i_1}\cap\ldots\cap W_{i_{\ell}}\right|\leq\Gamma_{\ell}$. Thus, $B_U$ is the event $\{N_n(U)>b\}$ and $\tilde B_U$ is the event $\{|N_n(U)|>b \mbox{ and, for every } \ell \in [k-1] \mbox{ and all distinct } i_1,...,i_\ell \in [k],\;
 |N_n(u_{i_1},...,u_{i_l})| \le \Gamma_\ell \}.$ In addition, note that
 $$
 \{X = 0\} = \bigwedge_U \neg B_U\,\,\text{ and }\,\,\{\tilde X = 0\} = \bigwedge_U\neg \tilde B_U.
 $$

Recall that $y$ is a constant. For a fixed $U$, after substituting
\begin{equation}
q=p^k,\quad a^+=b=np^k+x\sqrt{p^k(1-p^k)n\ln n},\quad x=\sqrt{2k}\left(1-\frac{\ln\ln n}{4k\ln n}-\frac{\ln\left[k!\sqrt{4\pi k}\right]-y}{2k\ln n}\right),
\label{eq:q_a_x}
\end{equation}
in~(\ref{central1}), we get
\begin{align}
{\sf P}(B_{U})={\sf P}(\xi_{n-k,p^k}>b)&\sim\frac{1}{2\sqrt{k\pi\ln n}}
\exp\left[-k\ln n\left(1-\frac{\ln\ln n}{4k\ln n}-\frac{\ln\left[k!\sqrt{4\pi k}\right]-y}{2k\ln n}\right)^2\right]\notag\\
&\sim\frac{1}{2\sqrt{k\pi\ln n}}
\exp\left[-k\ln n+\ln\sqrt{\ln n}+\ln\left[k!\sqrt{4\pi k}\right]-y\right]=\frac{k!}{n^k}e^{-y}.
\label{event_probability}
\end{align}


\subsection{The Janson-type inequality}
\label{Proof_Janson}

Denote
\begin{equation}\label{lambdas}
 \lambda=\sum_{U\in{[n]\choose k}}{\sf P}(B_U),\quad \tilde\lambda=\sum_{U\in{[n]\choose k}}{\sf P}(\tilde B_{U}),
\end{equation}
\begin{equation}
\Delta=\sum_{U_1,U_2\in{[n]\choose k}:\,U_1\cap U_2\neq\varnothing}{\sf P}(|N_n(U_1\cap U_2)|\leq\Gamma_{|U_1\cap U_2|},|N_n(U_1)|>b,|N_n(U_2)|>b).
\label{tildedelta}\end{equation}
\begin{lemma}
Under the conditions of Theorem~\ref{main}, the following bounds hold:
\begin{equation}
\mathrm{exp}\left[-\lambda+o(1)\right]\leq{\sf P}(X=0)\leq{\sf P}(\tilde X=0)\leq \mathrm{exp}\left[-(1+o(1))\tilde\lambda+(1+o(1))e^{\lambda}\Delta\right].
\label{main_bounds}
\end{equation}
\label{Janson-type}
\end{lemma}
{\it Proof}. The inequality ${\sf P}(X=0)\leq{\sf P}(\tilde X=0)$ follows from the definition of these random variables.

The inequality $e^{-\lambda+o(1)}\leq{\sf P}(X=0)$ follows from \cite[Theorem 6.3.3]{Alon} which is itself a direct corollary of the well known FKG inequality \cite[Theorem 6.2.1]{Alon} (note that, in the case $p=1/2$, this corollary is also known as Kleitman's Lemma~\cite{Kleitman}, see also \cite[Proposition 6.3.1]{Alon}). Thus, since the indicator random variables of the events $\neg B_{U}$ are decreasing functions of the edges of the random graph, we get
$$
 {\sf P}(X=0)={\sf P}\left(\bigwedge_{U}\neg B_{U}\right)\geq\prod_{U}(1-{\sf P}(B_{U}))=e^{\sum_{U}\ln(1-{\sf P}(B_{U}))}.
$$
From~(\ref{event_probability}), we get
\begin{equation}
{\sf P}(X=0)\geq e^{-\lambda+o(1)}.
\label{lower_exponent}
\end{equation}

The proof of the remaining inequality is close to the proof of Janson's inequality proposed by Boppona and Spencer~\cite{BopponaSpencer} (as well as to the proof of Suen's inequality proposed by Spencer~\cite{SpencerSuen}). However, it is harder since we need to overcome two difficulties. First, any two of $B$'s are not independent and, second, even after getting through the first barrier, we can not apply here the FKG inequality directly.

Let us consider an arbitrary ordering $\tilde B_1,\ldots,\tilde B_{n\choose k}$ of the events $\tilde B_U$. Then
\begin{equation}
 {\sf P}(\tilde X=0)={\sf P}\left(\bigwedge_{i=1}^{n\choose k} \neg \tilde B_i\right)=\prod_{i=1}^{n\choose k}\left[1-{\sf P}(\tilde B_i|\neg \tilde B_1\wedge\ldots\wedge\neg \tilde B_{i-1})\right].
\label{Janson_decomposition}
\end{equation}
Fix $i\in\big[{n\choose k}\big]$. Unfortunately, each event of $\tilde B_1,\ldots,\tilde B_{i-1}$ is not independent of $\tilde B_i$. Nevertheless, we may consider some slight modifications of them such that they are independent of $\tilde B_i$. Indeed, let $U$ define $\tilde B_i$ (i.e., $\tilde B_i=\tilde B_U$). Consider all the events (say, $\tilde B_1,\ldots,\tilde B_d$) among $\tilde B_1,\ldots,\tilde B_{i-1}$ (they are defined by $U_1,\ldots,U_d$ respectively) such that each of $U_1,\ldots,U_{d}$ has an empty intersection with $U$. Note that $d = d(i)$ depends on $i$. However, we write $d$ instead of $d(i)$ for shortening until the opposite is required. Let $j\in[d]$, and $\tilde B_j$ be defined by $U_j=\{u_1^j,\ldots, u_k^j\}$. In what follows, for $V\subset[n]$, we denote by $\hat N_n(V)$ the set of all common neighbors of $V$ in  $\mathcal{N}:=[n]\setminus U$.

For $\hat U=\{\hat u_1,\ldots,\hat u_k\}\in{\mathcal{N}\choose k}$, define
$$
 C_{\hat U}=\bigvee_{W_1,\ldots,W_k}\left\{\hat N_n(\hat u_1)=W_1,\ldots,\hat N_n(\hat u_k)=W_k\right\}.
$$
where the union is over sets $W_i\subset \mathcal{N}\setminus\{\hat u_i\}$, $i\in[k]$, such that $|W_1\cap\ldots\cap W_k|>b-k$ and, for every $\ell\in[k-1]$ and every distinct $i_1,\ldots,i_{\ell}\in[k]$, $\left|W_{i_1}\cap\ldots\cap W_{i_{\ell}}\right|\leq\Gamma_{\ell}$. Consider an arbitrary ordering $\hat U_1, \ldots, \hat U_{n-k\choose k}$ of the sets $\hat U \in {\mathcal{N}\choose k}$ such that $\hat U_j = U_j,$ $j\in[d].$  For $j\in\{1,\ldots,{n-k\choose k}\}$, set $C_j=C_{\hat U_j}$.
Clearly $C_1,\ldots,C_{n-k\choose k}$ do not depend on $\tilde B_i$ since $\tilde B_i$ is defined only by edges having vertices in $U$. Note also that, generally speaking, the events $C_j$ depend on $i$ by the definition of $\mathcal{N},$ but, hereinafter, we write $C_j$ instead of $C_j(i)$ to prevent overloading with double indexations.

By the proven lower bound~(\ref{lower_exponent}) (note that this lower bound still holds true if we replace $b$ with $b+O(1)$ in the definition of $X$),
\begin{align}
 {\sf P}\left(\bigwedge_{i=1}^{{n-k\choose k}}\neg C_i\right)&\geq
 {\sf P}\left(\text{there is no }\hat U\in{\mathcal{N}\choose k}\text{ in }G(n-k,p)\text{ with }|\hat N_n(\hat U)|>b-k\right)\notag\\
 &\geq {\sf P}\left(\text{there is no }\hat U\in{[n]\choose k}\text{ in }G(n,p)\text{ with }|N_n(\hat U)|>b-k\right) \geq e^{-\lambda +o(1)}.
\label{auxiliary_events}
\end{align}
Clearly,
\begin{equation}
 {\sf P}(\tilde B_i\wedge\neg C_1\wedge\ldots\wedge\neg C_d)\leq{\sf P}(\tilde B_i\wedge\neg \tilde B_1\wedge\ldots\wedge\neg \tilde B_d)
\label{BB_to_BC}
\end{equation}
and
$$
{\sf P}(\neg C_1\wedge\ldots\wedge\neg C_d)\geq {\sf P}(\neg \tilde B_1\wedge\ldots\wedge\neg \tilde B_d)-{\sf P}\left(\exists V\in{[n]\choose k}\text{ s.t. }|N_n(V)|\in\{b-k+1,\ldots,b\}\right)-
$$
\begin{equation}
\sum_{\ell=1}^{k-1}{\sf P}\left(\exists V\in{[n]\choose\ell}\text{ s.t. }|N_n(V)|\in\{\Gamma_{\ell}+1,\ldots,\Gamma_{\ell}+k\}\right)={\sf P}(\neg\tilde B_1\wedge\ldots\wedge\neg\tilde B_d)-o\left(\frac{1}{\ln n}\right)
\label{B_to_C}
\end{equation}
by the union bound (the remainder $o(1/\ln n)$ does not depend on $i$). Indeed, for every $a\in\{b-k+1,\ldots,b\}$ and every $V\in{[n]\choose k}$, from~(\ref{demoivre}), we get (here, we apply (\ref{demoivre}) with the parameters $q,a,x$ defined in (\ref{eq:q_a_x}); note that additive $O(1)$ terms that come with $n$ and $a$ here do not affect asymptotics)
$$
{\sf P}\left(|N_n(V)|=a\right)={\sf P}\left(\xi_{n-k,p^k}=a\right)\sim
 \frac{\sqrt{2k\ln n}k!e^{-y}}{n^k\sqrt{np^k(1-p^k)}}= o\left(\frac{1}{n^k \ln n}\right).
$$
For every $\ell\in\{1,\ldots,k-1\}$, $V\in{[n]\choose\ell}$, $a\in\{\Gamma_{\ell}+1,\ldots,\Gamma_{\ell}+k\}$, after substituting
$$
q=p^{\ell},\quad a=np^{\ell}+x\sqrt{p^{\ell}(1-p^{\ell})n\ln n},\quad x=\sqrt{2\ell}+o\left(n^{-1/(2k)}\right)
$$
in~(\ref{demoivre}), we get
$$
{\sf P}\left(|N_n(V)|=a\right)={\sf P}\left(\xi_{n-\ell,p^{\ell}}=a\right)\sim\frac{1}{\sqrt{2\pi n p^{\ell}(1-p^{\ell})}}e^{-\ell \ln n}=o\left(\frac{1}{n^{\ell}(\ln n)^{3/2} }\right).
$$

Below, we use a standard tool from the proof of Janson's inequality:
$$
{\sf P}(\tilde B_i|\neg \tilde B_1\wedge\ldots\wedge\neg \tilde B_{i-1})\geq
\frac{{\sf P}(\tilde B_i\wedge\neg \tilde B_1\wedge\ldots\wedge\neg \tilde B_d\wedge\neg \tilde B_{d+1}\wedge\ldots\wedge\neg \tilde B_{i-1})}{{\sf P}(\neg \tilde B_1\wedge\ldots\wedge\neg \tilde B_d)}=
$$
$$
 {\sf P}(\tilde B_i|\neg \tilde B_1\wedge\ldots\wedge\neg \tilde B_d)
 {\sf P}(\neg \tilde B_{d+1}\wedge\ldots\wedge\neg \tilde B_{i-1}|\tilde B_i\wedge\neg \tilde B_1\wedge\ldots\wedge\neg \tilde B_d)=
 $$
$$
 \frac{{\sf P}(\neg C_1\wedge\ldots\wedge\neg C_d)}{{\sf P}(\neg B_1\wedge\ldots\wedge\neg \tilde B_d)}\times\frac{{\sf P}(\tilde B_i\wedge\neg \tilde B_1\wedge\ldots\wedge\neg \tilde B_d)}{{\sf P}(\tilde B_i\wedge\neg C_1\wedge\ldots\wedge\neg C_d)}\times{\sf P}(\tilde B_i)\times
$$
\begin{equation}
 {\sf P}(\neg \tilde B_{d+1}\wedge\ldots\wedge\neg \tilde B_{i-1}|\tilde B_i\wedge\neg \tilde B_1\wedge\ldots\wedge\neg \tilde B_d)
\label{std_tool}
\end{equation}
since $\tilde B_i$ is independent of $\neg C_1\wedge\ldots\wedge\neg C_d$. It remains to estimate the factor from the last line.  For every $j\in\{d+1,\ldots,i-1\}$, from~(\ref{auxiliary_events})~and~(\ref{BB_to_BC}),
$$
 {\sf P}(\tilde B_j|\tilde B_i\wedge\neg\tilde B_1\wedge\ldots\wedge\neg\tilde B_d)
 =\frac{{\sf P}(\tilde B_j\wedge\tilde B_i\wedge\neg \tilde B_1\wedge\ldots\wedge\neg\tilde  B_d)}{{\sf P}(\tilde B_i\wedge\neg\tilde B_1\wedge\ldots\wedge\neg\tilde B_d)}
 \leq\frac{{\sf P}(\tilde B_j\wedge\tilde B_i\wedge\neg \tilde B_1\wedge\ldots\wedge\neg\tilde  B_d)}{{\sf P}(\tilde B_i){\sf P}(\neg C_1\wedge\ldots\wedge\neg C_d)}\leq
$$
\begin{equation*}
  \frac{{\sf P}(\tilde B_j\wedge\tilde B_i)}{{\sf P}(\tilde B_i){\sf P}(\neg C_1\wedge\ldots\wedge\neg C_{{n-k\choose k}})}\leq {\sf P}(\tilde B_j|\tilde B_i)e^{\lambda+o(1)}.
 \label{Janson_help}
\end{equation*}
Therefore,
$$
 {\sf P}(\neg \tilde B_{d+1}\wedge\ldots\wedge\neg \tilde B_{i-1}|\tilde B_i\wedge\neg \tilde B_1\wedge\ldots\wedge\neg \tilde B_d)\geq
$$
\begin{equation}
 1-\sum_{j=d+1}^{i-1}{\sf P}(\tilde B_j|\tilde B_i\wedge\neg\tilde B_1\wedge\ldots\wedge\neg\tilde B_d)\geq
 1-e^{\lambda+o(1)}\sum_{j=d+1}^{i-1}{\sf P}(\tilde B_j|\tilde B_i).
\label{Janson_help1}\end{equation}
Estimating from below the first, second and fourth multipliers in the right-hand side of (\ref{std_tool}) by (\ref{B_to_C}), (\ref{BB_to_BC}) and (\ref{Janson_help1}), respectively, we derive the following lower bound
$$
{\sf P}(\tilde B_i|\neg \tilde B_1\wedge\ldots\wedge\neg \tilde B_{i-1}) \geq \left(1-\frac{o(1/\ln n)}{{\sf P}(\neg \tilde B_1\wedge\ldots\wedge\neg \tilde B_{n\choose k})}\right){\sf P}(\tilde B_i)\left[1-e^{\lambda+o(1)}\sum_{j=d(i)+1}^{i-1}{\sf P}(\tilde B_j|\tilde B_i)\right].
 $$
Thus, combining the latter and (\ref{Janson_decomposition}) and recalling that $\lambda = O(1),$ we finally get
$$
 {\sf P}(\tilde X=0)\leq\prod_{i=1}^{n\choose k}\left(1-\left(1-\frac{o(1/\ln n)}{{\sf P}(\neg \tilde B_1\wedge\ldots\wedge\neg \tilde B_{n\choose k})}\right){\sf P}(\tilde B_i)\left[1-e^{\lambda+o(1)}\sum_{j=d(i)+1}^{i-1}{\sf P}(\tilde B_j|\tilde B_i)\right]\right)\leq
$$
$$
 \mathrm{exp}\left[-\left(1-\frac{o(1/\ln n)}{\mathrm{exp}[-\lambda+o(1)]}\right)\left(\sum_{i=1}^{n\choose k}{\sf P}(\tilde B_i)-e^{\lambda+o(1)}\sum_{i=1}^{n\choose k}\sum_{j=d(i)+1}^{i-1}{\sf P}(\tilde B_j\wedge \tilde B_i)\right)\right]\leq
 $$
 $$
 \mathrm{exp}\left[-(1+o(1))\left(\tilde\lambda-e^{\lambda+o(1)}\Delta\right)\right].
$$
A noteworthy detail is that $\sum_{i=1}^{n\choose k}\sum_{j=d(i)+1}^{i-1}{\sf P}(\tilde B_j\wedge\tilde B_i)$ is strictly less than $\Delta$, since, in the definition of $\Delta$, we remove the restrictions on the cardinalities of sets of common neighbors of all proper subsets of both $k$-sets but the only subset which is the intersection of $k$-sets. So, the upper bound in~(\ref{main_bounds}) can be strengthened, but, for our purpose, this bound is more convenient. $\Box$

\subsection{The second moment}
\label{Proof_Moments}

In this section, we prove the following
\begin{lemma}
Let $\lambda,$ $\tilde \lambda$ and $\Delta$ be defined by (\ref{lambdas}) and (\ref{tildedelta}). Then under the conditions of Theorem 1,
$$
\tilde\lambda\sim\lambda\sim e^{-y}, \text{ and }\Delta\to 0
$$
as $n\to\infty.$
\label{parameters_estimation}
\end{lemma}

{\it Proof.} Let us start from estimating $\lambda,\tilde\lambda$ and showing that $\tilde\lambda=\lambda(1+o(1))$. From~(\ref{event_probability}),
\begin{equation}
\lambda =\sum_{U\in{[n]\choose k}}{\sf P}(B_U)={n\choose k}\frac{k!}{n^k}e^{-y}(1+o(1))=e^{-y}(1+o(1)).
\label{expectation}
\end{equation}


\subsubsection{Estimation of $\tilde\lambda$}
\label{lambda}

Fix a $k$-set $U =\{u_1, \ldots, u_k\}\in{[n]\choose k}$. 
$$
{\sf P}(B_U) \geq {\sf P}(\tilde B_U)\geq
{\sf P}(B_U) - \sum_{\ell=1}^{k-1}\sum_{V\in{U\choose\ell}}{\sf P}\left(|N_n(U)| > b,\,|N_n(V)|>\Gamma_{\ell}\right)=
$$
\begin{equation}
{\sf P}(B_U) - \sum_{\ell=1}^{k-1}{k\choose\ell}{\sf P}\left(|N_n(U)| > b,\,|N_n(\{u_1,\ldots,u_{\ell}\})|>\Gamma_{\ell}\right).
\label{tilde_B_estimation}
\end{equation}
Denote $U_{\ell}=\{u_1,\ldots,u_{\ell}\}$.
Thus, we get
\begin{align*}
{\sf P}(|N_n(U)| > b,\,|N_n(U_{\ell})|>\Gamma_{\ell})
& \leq \sum_{i>\Gamma_{\ell}}{\sf P}(\xi_{n-\ell,p^{\ell}}=i){\sf P}(\xi_{i, p^{k-\ell}}>b-(k-\ell))\\
&\leq \sum_{\Gamma_{\ell}<i\leq np^{\ell}+\sqrt{2kp^{\ell}(1-p^{\ell})n\ln n}}{\sf P}(\xi_{n-\ell,p^{\ell}}=i){\sf P}(\xi_{i, p^{k-\ell}}>b-(k-\ell))\\
&+{\sf P}\left(\xi_{n-\ell,p^{\ell}}>np^{\ell}+\sqrt{2kp^{\ell}(1-p^{\ell})n\ln n}\right),
\end{align*}
We aim at proving that the latter sum is $o(n^{-k}).$
From relations (\ref{demoivre}) and (\ref{central1}) (since $x$ in this case belongs to the interval $\left(\sqrt{2\ell}, \sqrt{2k}\right]$, the conditions (\ref{assumptions_moivre}) and (\ref{assumptions_moivre1}) immediately follow from \eqref{cond}, as usual),
\begin{multline} \label{additional_and_similar}
 {\sf P}(\xi_{n-\ell,p^{\ell}}=i)=\frac{\mathrm{exp}\left[-\frac{(np^{\ell}-i)^2}{2np^{\ell}(1-p^{\ell})}\right]}{\sqrt{2\pi np^{\ell}(1-p^{\ell})}}(1+o(1))\\\quad\text{uniformly over }i\in\left(\Gamma_{\ell},np^{\ell}+\sqrt{2kp^{\ell}(1-p^{\ell})n\ln n}\right],
\end{multline}
\begin{equation} \label{additional_and_similar2}
 {\sf P}\left(\xi_{n-\ell,p^{\ell}}>np^{\ell}+\sqrt{2kp^{\ell}(1-p^{\ell})n\ln n}\right)=
\frac{1}{2\sqrt{k\pi\ln n}}n^{-k}(1+o(1)).
\end{equation}

It remains to estimate ${\sf P}(\xi_{i, p^{k-\ell}}>b-(k-\ell))$. 
Let us verify the conditions of~(\ref{central1}). Since $k-\ell$ is constant, it is sufficient to prove that
\begin{equation}
1\ll b-ip^{k-\ell}=\Theta\left(\sqrt{\frac{\ln n}{\ln i}}\right)\times\sqrt{p^{k-\ell}(1-p^{k-\ell})i\ln i }.
\label{formula:referee}
\end{equation}
Indeed, we should check \eqref{assumptions_moivre} and \eqref{assumptions_moivre1} for
$$
n = i,\quad q = p^{k-\ell}\,\,\text{ and }\,\,x = (b-ip^{k-\ell})/\sqrt{p^{k-\ell}(1-p^{k-\ell})i\ln i }.
$$
From (\ref{formula:referee}) it follows that $0<x = \Theta(\sqrt{\ln n/\ln i})=\Theta(1)$ which immediately implies the first two conditions in \eqref{assumptions_moivre1}. The rest is straightforward due to the restrictions on $p$ given in Theorem~\ref{main} and the relation $x=\Theta(1)$. 

Now, let us prove \eqref{formula:referee}. For $i$ in the range,
$$
 b-ip^{k-\ell}\geq\sqrt{2k p^k(1-p^k)n\ln n}(1+o(1))-\sqrt{2kp^{2k-\ell}(1-p^{\ell})n\ln n}=
$$
\begin{equation}
 \sqrt{2kp^kn\ln n}\left(\sqrt{1-p^k}(1+o(1))-\sqrt{p^{k-\ell}-p^{k}}\right)\geq\sqrt{2kp^k n\ln n}\sqrt{1-p^k}\left(1-\sqrt{\ell/k}+o(1)\right),
\label{b-i_lower}
\end{equation}
where the latter inequality follows from
\begin{equation}
\sqrt{\frac{k(p^{k-\ell}-p^k)}{1-p^k}}<\sqrt{\ell}.
\label{auxiliary_inequality}
\end{equation}
Indeed,
$$
\frac{\partial}{\partial p}\left[\ell(1-p^k)-k(p^{k-\ell}-p^k)\right]=
k(\ell-k)p^{k-\ell-1}(1-p^{\ell})<0
$$
for $p\in(0,1)$, and $\left.\ell(1-p^k)-k(p^{k-\ell}-p^k)\right|_{p=1}=0$.

From (\ref{b-i_lower}), it follows that $b-ip^{k-\ell}\gg 1$. Moreover,
$$
 \frac{b-ip^{k-\ell}}{\sqrt{p^{k-\ell}(1-p^{k-\ell})i\ln i}}\geq \sqrt{\frac{2k(1-p^k)\ln n}{(1-p^{k-\ell})\ln i}}\left(1-\sqrt{\ell/k}+o(1)\right)\geq
 \sqrt{\frac{\ln n}{\ln i}}\left(\sqrt{2k}-\sqrt{2\ell}+o(1)\right)
$$
as needed.

Let us verify the upper bound:
$$
 \frac{b-ip^{k-\ell}}{\sqrt{p^{k-\ell}(1-p^{k-\ell})i\ln i}}\leq\frac{\sqrt{2kp^k(1-p^k)n\ln n}}{\sqrt{p^{k}(1-p^{k-\ell})n\ln i}(1+o(1))}=O\left(\sqrt{\frac{\ln n}{\ln i}}\right),
$$ that completes the proof of \eqref{formula:referee}.

Therefore, we may apply~(\ref{central1}) for ${\sf P}(\xi_{i, p^{k-\ell}}>b-(k-\ell))$ as well: 
\begin{equation} \label{bkell}
{\sf P}(\xi_{i, p^{k-\ell}}>b-(k-\ell))=\frac{\sqrt{ip^{k-\ell}(1-p^{k-\ell})}}{\sqrt{2\pi}(b-ip^{k-\ell})}
\exp\left[-\frac{(b-ip^{k-\ell})^2}{2ip^{k-\ell}(1-p^{k-\ell})}\right](1+o(1))\leq
\end{equation}
$$
\frac{\sqrt{1-p^{k-\ell}}}{2\sqrt{\pi k\ln n}\left(\sqrt{1-p^k}-\sqrt{p^{k-\ell}-p^{k}}\right)}
\exp\left[-\frac{(b-ip^{k-\ell})^2}{2ip^{k-\ell}(1-p^{k-\ell})}\right](1+o(1)).
$$

Putting it all together, we get
$$
{\sf P}(|N_n(U)| > b,\,|N_n(U_{\ell})|>\Gamma_{\ell}) \leq
$$
\begin{equation}
\left(\frac{\sqrt{1-p^{k-\ell}}}{2\pi\sqrt{2k p^{\ell}(1-p^{\ell})n\ln n}\left(\sqrt{1-p^k}-\sqrt{p^{k-\ell}-p^k}\right)}\sum e^{-f(i)}+
\frac{1}{2\sqrt{k\pi\ln n}}n^{-k}\right)(1+o(1)),
\label{lambda_reminder}
\end{equation}
where the summation is over $i\in\left(\Gamma_{\ell},np^{\ell}+\sqrt{2kp^{\ell}(1-p^{\ell})n\ln n}\right]$ and
$$
f(i)=\frac{\left(np^{\ell}-i\right)^2}{2np^{\ell}\left(1-p^{\ell}\right)}+\frac{\left(b-ip^{k-\ell}\right)^2}{2ip^{k-\ell}\left(1-p^{k-\ell}\right)}.
$$
Let $i>\Gamma_{\ell}$. Denote $i =np^{\ell} + x\sqrt{np^{\ell}(1-p^{\ell})\ln n}.$ Then, omitting direct, but tedious calculations, we have

$$
f(i)=\frac{x^2(1-p^k) -2\sqrt{2k}x\sqrt{(p^{k-\ell} - p^k)(1-p^k)} + 2k(1-p^k)}{2(1-p^{k-\ell})}\ln n +
$$
 \begin{equation}
 \frac{x \sqrt{(p^{k-\ell} - p^k)(1 - p^k)} - \sqrt{2k} (1 - p^k)}{2\sqrt{2k} (1 - p^{k-\ell})}\ln\ln n (1 + o(1)),
\label{f(i)}
\end{equation}
From the definition of $\Gamma_{\ell}$, $x>\sqrt{2\ell}$. On the other hand, the function of $x$ in the numerator of the first summand in the right-hand side of~\eqref{f(i)} (denote it by $\hat{f}(x)$) achieves its minimum at
$\sqrt{\frac{2k(p^{k-\ell}-p^k)}{1-p^k}}<\sqrt{2\ell}$ due to (\ref{auxiliary_inequality}). Therefore, $\hat{f}(x)>\hat{f}(\sqrt{2\ell}).$

Next, let us prove that the ratio before $\ln\ln n$ in the right-hand side of~\eqref{f(i)} is negative and bounded. Since $x\le \sqrt{2k}$, we have
$$
\frac{x \sqrt{(p^{k-\ell} - p^k)(1 - p^k)} - \sqrt{2k} (1 - p^k)}{2\sqrt{2k} (1 - p^{k-\ell})} \le \frac{\sqrt{2k (1 - p^k)} \big(\sqrt{p^{k-\ell} - p^k} - \sqrt{1 - p^k}\big)}{2\sqrt{2k} (1 - p^{k-\ell})} < 0.
$$
On the other hand,
$$
\frac{x \sqrt{(p^{k-\ell} - p^k)(1 - p^k)} - \sqrt{2k} (1 - p^k)}{2\sqrt{2k} (1 - p^{k-\ell})} > -\frac{(1 - p^k)}{2(1 - p^{k-\ell})} > - \frac{\frac{k}{k-\ell} (1 - p^{k-\ell})}{2(1 - p^{k-\ell})} = - \frac{k}{2(k-\ell)},
$$ where the second inequality follows from the relation $(k-\ell)(1-p^k)<k(1 - p^{k-\ell})$ which can be proved similarly to (\ref{auxiliary_inequality}).

To prove \eqref{33}, we need to make the bound (\ref{auxiliary_inequality}) tighter. Namely, we want to show that
\begin{equation} \label{aux2}
k\frac{p^{k-\ell} - p^k}{1 - p^k} < \ell p^{(k-\ell)/2}.
\end{equation}
Indeed, it is equivalent to show that
$$
h(p) := \ell(1 - p^k) - k(p^{(k-\ell)/2} - p^{(k+\ell)/2})>0.
$$
Since $h(1) = 0,$ it is sufficient to prove that $h^\prime(p) <0$ for $p\in(0,1).$ We have,
$$
h^\prime(p) = -\ell k p^{k-1} - k (k-\ell)p^{(k-\ell)/2-1}/2 + k (k+\ell) p^{(k+\ell)/2-1}/2 < 0,
$$
or, equivalently,
\begin{equation} \label{aux3} 2 \ell p^{(k+\ell)/2} + (k-\ell) - (k+\ell) p^\ell > 0.\end{equation}
Its derivative $\ell (k+\ell) p^{(k+\ell)/2-1} - \ell (k+\ell) p^{\ell-1} < 0$ for $p\in(0,1)$ and $2 \ell p^{(k+\ell)/2} + (k-\ell) - (k+\ell) p^\ell|_{p=1} = 0,$ thus \eqref{aux3} and, consequently, \eqref{aux2} hold.

 Thus, after plugging in the numerator of the first summand in the right-hand side of~(\ref{f(i)}) $\sqrt{2\ell}$ and using the relation \eqref{aux2}, we get that
\begin{eqnarray}
f(i)=\left[k + \frac{(1-p^k)\left(x - \frac{\sqrt{2k(p^{k-\ell}-p^k)}}{\sqrt{1-p^k}}\right)^2}{2(1-p^{k-\ell})}\right]\ln n-O(\ln\ln n)\nonumber\\
>
\left[k + \frac{1}{2}\left(x - \frac{\sqrt{2k(p^{k-\ell}-p^k)}}{\sqrt{1-p^k}}\right)^2\right]\ln n-O(\ln\ln n)\nonumber\\
>\left[k + \left(\sqrt{2\ell} - \sqrt{2\ell p^{(k-\ell)/2}}\right)^2/2\right]\ln n-O(\ln\ln n)\nonumber\\
=\left[k + \ell (1 - p^{(k-\ell)/4})^2\right]\ln n-O(\ln\ln n).\label{33}
\end{eqnarray}
Using the assumption $1-p\gg\sqrt{\frac{\ln\ln n}{\ln n}}$, we have $(1 - p^{(k-\ell)/4})^2 \gg \frac{\ln\ln n}{\ln n}$. Thus,
$$
f(i) > \ln n (k + \omega(\ln\ln n/\ln n)) - O(\ln\ln n) = k\ln n + \omega(\ln\ln n).
$$

Thus,
\begin{equation}
\sum_{i\in\left(\Gamma_{\ell},np^{\ell}+\sqrt{2kp^{\ell}(1-p^{\ell})n\ln n}\right]} e^{-f(i)}<
\sqrt{2kp^{\ell}(1-p^{\ell})n\ln n}n^{-k}.
\label{a_little_problem}
\end{equation}
From~(\ref{tilde_B_estimation}),~(\ref{lambda_reminder}) and (\ref{a_little_problem}), we get
$$
{\sf P}(B_U) \geq {\sf P}\left(\tilde B_U\right)\geq{\sf P}(B_U) - o(n^{-k}).
$$
Finally, from~(\ref{event_probability})~and~(\ref{expectation}),
$$
 e^{-y}+o(1)={n\choose k}\left({\sf P}(B_U)-o(n^{-k})\right)\leq
{n\choose k}{\sf P}\left(\tilde B_U\right)=\tilde\lambda\leq\lambda={n\choose k}{\sf P}(B_U)=e^{-y}+o(1).
$$

\subsubsection{Estimation of $\Delta$}
\label{delta}

It remains to estimate from above $\Delta$ and prove its convergence to $0$. By the definition,
$$
\Delta = \sum_{\ell=1}^{k-1}\sum_{V\in{[n]\choose\ell},\,U_1,U_2\in{[n]\choose k}:U_1\cap U_2=V}{\sf P}(|N_n(V)|\leq\Gamma_{\ell},|N_n(U_1)|>b,|N_n(U_2)|>b).
$$
Fix $k$-sets $U_1,U_2$ such that $|U_1\cap U_2|=\ell$, $1\leq\ell\leq k-1$. Set $V=U_1\cap U_2$. Let
$$
 A^{\ell}_{U_1,U_2}=\{|N_n(V)|\leq\Gamma_{\ell},|N_n(U_1)|>b,|N_n(U_2)|>b\}.
$$
Then
\begin{equation}
\Delta = \sum_{\ell=1}^{k-1}{n\choose\ell}{n-\ell\choose k-\ell}{n-k\choose k-\ell}{\sf P}(A^{\ell}_{U_1,U_2})=\sum_{\ell=1}^{k-1}O\left(n^{2k-\ell}{\sf P}(A^{\ell}_{U_1,U_2})\right).
\label{delta_A}
\end{equation}
So it is sufficient to show that ${\sf P}(A^{\ell}_{U_1,U_2}) =o(n^{-(2k-\ell)})$. Obviously,
$$
{\sf P}(A^{\ell}_{U_1,U_2}) \leq \sum_{i\leq\Gamma_{\ell}} {\sf P}(\xi_{n, p^{\ell}} = i) [{\sf P}(\xi_{i, p^{k-\ell}}>b-(k-\ell))]^2 \leq
$$
$$
\sum_{np^{\ell}-2\sqrt{k p^{\ell}(1-p^{\ell})n\ln n}< i\leq\Gamma_{\ell}} {\sf P}(\xi_{n, p^{\ell}} = i) [{\sf P}(\xi_{i, p^{k-\ell}}>b-(k-\ell))]^2 + $$
$$
{\sf P}\left(\xi_{n,p^{\ell}}\leq np^{\ell}-2\sqrt{kp^{\ell}(1-p^{\ell})n\ln n}\right).
$$
From relations (\ref{demoivre}) and ($10'$), we get
$$
{\sf P}(\xi_{n,p^{\ell}}=i)=\frac{\exp\left[-\frac{(np^{\ell}-i)^2}{2np^{\ell}(1-p^{\ell})}\right]}{\sqrt{2\pi np^{\ell}(1-p^{\ell})}}(1+o(1))\text{ uniformly over }i\in\left(np^{\ell}-2\sqrt{kp^{\ell}(1-p^{\ell})n\ln n}, \Gamma_{\ell}\right],
$$
$$
{\sf P}\left(\xi_{n,p^{\ell}}\leq np^{\ell}-2\sqrt{kp^{\ell}(1-p^{\ell})n\ln n}\right) = \frac{1}{2\sqrt{2k\pi \ln n}} n^{-2k}(1 + o(1)) = o(n^{-2k}).
$$

As in Section~\ref{lambda}, to estimate ${\sf P}(\xi_{i, p^{k-\ell}}>b-(k-\ell))$, we should verify the conditions of (\ref{central1}), namely, the assumptions (\ref{assumptions_moivre}) and (\ref{assumptions_moivre1}). For $i\in\left(np^{\ell}-2\sqrt{kp^{\ell}(1-p^{\ell})n\ln n}, \Gamma_{\ell}\right]$, the value of $b-ip^{k-\ell}$ is even bigger than in (\ref{b-i_lower}) and, therefore, here, (\ref{b-i_lower}) holds as well. Moreover,
$$
 \frac{b-ip^{k-\ell}}{\sqrt{p^{k-\ell}(1-p^{k-\ell})i\ln i }}\leq\frac{\sqrt{2kp^k(1-p^k)n\ln n}+2\sqrt{kp^k(p^{k-\ell}-p^k)n\ln n}}{\sqrt{p^{k}(1-p^{k-\ell})n\ln i }(1+o(1))}=O\left(\sqrt{\frac{\ln n}{\ln i}}\right).
$$
Therefore, similarly to corresponding argument in Section~\ref{lambda}, we may apply~(\ref{central1}) and get, in the same way as in~(\ref{bkell}),
$$
{\sf P}(\xi_{i, p^{k-\ell}}>b-(k-\ell)) \leq
\frac{\sqrt{1-p^{k-\ell}}}{\sqrt{2\pi\ln n}\left(\sqrt{2k(1-p^k)}-\sqrt{2\ell(p^{k-\ell}-p^{k})}\right)}e^{-\frac{(b-ip^{k-\ell})^2}{2ip^{k-\ell}(1-p^{k-\ell})}}(1+o(1)).
$$
Therefore, we get from above
$$
{\sf P}(A_{U_1,U_2}^{\ell})\leq
$$
\begin{equation}
(1+o(1))\frac{1-p^{k-\ell}}{(2\pi)^{3/2}\ln n\sqrt{n p^\ell(1-p^{\ell})}\left(\sqrt{2k(1-p^k)}- \sqrt{2\ell p^{k-\ell}(1-p^{\ell})}\right)^2}\sum e^{-g(i)}+ o(n^{-2k}),
\label{delta_A_estimation}
\end{equation}
where the summation is over $i\in\left(np^{\ell}-2\sqrt{kp^{\ell}(1-p^{\ell})n\ln n}, \Gamma_{\ell}\right]$ and
$$
g(i)=\frac{(np^{\ell}-i)^2}{2np^{\ell}(1-p^{\ell})}+\frac{(i p^{k-\ell}-b)^2}{i p^{k-\ell}(1-p^{k-\ell})}.
$$

As in Section~\ref{lambda}, denote $i = np^{\ell} + x\sqrt{np^{\ell}(1-p^{\ell})\ln n}.$ Notice that $x \in (-2\sqrt{k}, \sqrt{2\ell}].$

We get
$$
g(i) =
\frac{x^2(1+p^{k-\ell} - 2p^k) - 4\sqrt{2k}\sqrt{(p^{k-\ell}-p^k)(1-p^k)}x + 4k(1-p^k)}{2(1-p^{k-\ell})}\ln n +
$$
$$
\frac{x \sqrt{(p^{k-l} - p^k)(1 - p^k)} - \sqrt{2k} (1 - p^k)}{\sqrt{2k} (1 - p^{k-l})}\ln\ln n (1 + o(1)) =: \tilde g_p(x) \ln n + \hat g_p(x) \ln\ln n (1+o(1)).
$$
We note that $\hat g_p(x)/2$ is equal to the ratio before $\ln\ln n$ in the right-hand side of (\ref{f(i)}). 
It can be proved similarly as in Section~\ref{lambda} that $\hat g_p(x)$ is negative and bounded from below.

As in Section~\ref{lambda}, in order to show ${\sf P}(A^{\ell}_{U_1,U_2}) =o(n^{-(2k-\ell)}),$ we prove $\tilde g_p(x) \ge 2k - \ell + \omega(\ln\ln n/\ln n).$ The minimum of $\tilde g_p(x)$ equals $\tilde g_p(x_0) = \frac{2k(1-p^k)}{1+p^{k-\ell} - 2p^k}$, where
$$
x_0 = \frac{2\sqrt{2k}\sqrt{(p^{k-\ell}-p^k)(1-p^k)}}{1+p^{k-\ell} - 2p^k}.
$$
Note that $\tilde g_p(x_0)$ decreases in $p\in(0,1)$ since \[\frac{2k}{\tilde g_p(x_0)} = 2 - \frac{1 - p^{k-\ell}}{1 - p^k},\] and $(1 - p^{k-\ell})/(1 - p^k)$ decreases in $p\in(0,1)$ as well.
Let $\varepsilon\in\big(0,\frac{2}{k-1}\big)$ and $p_0\in(0,1)$ be the root of the equation $\frac{1-p^k}{p^{k-\ell}-p^k} = \frac{2k-\ell}{\ell} + \varepsilon$. Then we have 
$$
\tilde g_{p_0}(x_0) = 2k \left(\frac{p_0^{k-\ell} - p_0^k}{1 - p_0^k} + 1\right)^{-1} = 2k\left(\frac{\ell}{2k - \ell(1 - \varepsilon)} + 1\right)^{-1} = 2k\frac{2k - (1-\varepsilon) \ell}{2k + \varepsilon \ell}.
$$
Thus, $$\tilde g_p(x_0) > 2k\frac{2k - (1-\varepsilon) \ell}{2k + \varepsilon \ell} = 2k-\ell + \frac{\varepsilon \ell^2}{2k + \varepsilon \ell}> 2k-\ell$$ for all $p<p_0.$ For such $p$, we immediately get the required relation ${\sf P}(A^{\ell}_{U_1,U_2}) =o(n^{-(2k-\ell)})$.

If $p\geq p_0$, then $\frac{1-p^k}{p^{k-\ell}-p^k}\leq\frac{2k-\ell}{\ell} + \varepsilon,$ since the function $\frac{1-p^k}{p^{k-\ell}-p^k}$ decreases for $p\in(0,1).$ Therefore, since the function $x_0 = x_0(p)$ is increasing for $p\in (0,1),$ we get
\begin{equation} \label{x_0}
x_0 = \frac{2\sqrt{2k}\sqrt{(p^{k-\ell}-p^k)(1-p^k)}}{1+p^{k-\ell} - 2p^k} \geq x_0(p_0) = 2\sqrt{2k} \frac{\sqrt{\ell(2k - (1-\varepsilon)\ell)}}{2k + \varepsilon\ell} >\sqrt{2\ell},
\end{equation}
where the last inequality follows from the relation $4k(2k - (1-\varepsilon)\ell) > (2k + \varepsilon\ell)^2$ which can be shown directly using $\ell \le k-1$ and $\varepsilon < 2/(k-1).$
Moreover,
$$
\tilde g_p\left(\sqrt{2\ell}\right) = \frac{\ell(1 + p^{k-\ell} - 2p^k) - 2\sqrt{2k}\sqrt{2\ell}\sqrt{(p^{k-\ell}-p^k)(1-p^k)} + 2k(1-p^k)}{1-p^{k-\ell}} =
$$
$$
\frac{\left[\sqrt{2\ell(1-p^k)} - \sqrt{2k(p^{k-\ell} - p^k)}\right]^2 - \ell(1-p^{k-\ell}) + 2k(1-p^{k-\ell})}{1-p^{k-\ell}} = 2k - \ell +\omega\left(\frac{\ln\ln n}{\ln n}\right)
$$
since $\left[\sqrt{2\ell(1-p^k)} - \sqrt{2k(p^{k-\ell} - p^k)}\right]^2/(1-p^{k-\ell})=\omega\left(\frac{\ln\ln n}{\ln n}\right)$ for $1-p\gg\sqrt{\frac{\ln\ln n}{\ln n}},$ see \eqref{33} and the argument below \eqref{33}. As $x\leq\sqrt{2\ell}< x_0$, we get $\tilde g_p(x)\geq \tilde g_p\left(\sqrt{2\ell}\right).$ Since $\hat g_p(x)=O(1)$, 
\begin{equation}
\sum_{i\in\left(np^{\ell}-2\sqrt{kp^{\ell}(1-p^{\ell})n\ln n}, \Gamma_{\ell}\right]} e^{-g(i)}<4\sqrt{kp^{\ell}(1-p^{\ell})n\ln n}n^{-(2k-\ell)}.
\label{delta_small_sum}
\end{equation}
Relations~(\ref{delta_A}),~(\ref{delta_A_estimation}),~(\ref{delta_small_sum}) imply the desired convergence $\Delta\to 0$. $\Box$

\subsection{Final steps}

From Lemma~\ref{Janson-type} and Lemma~\ref{parameters_estimation} we immediately get the statement of Theorem~\ref{main} for $m=1$ since ${\sf P}(\Delta^1_{k;n}\leq b)={\sf P}(X_n=0)$.

The very last step is to find the limit of ${\sf P}(\Delta^m_{k;n}\leq b)-{\sf P}(\Delta^{m-1}_{k;n}\leq b)={\sf P}(X_n=m-1)$ for $m\geq 2$.

Note that the probability $\delta_n$ of the existence of two distinct overlapping $k$-sets $U_1,U_2\in{[n]\choose k}$ such that $|N_n(U_j)|>b$ for both $j=1$ and $j=2$ is at most $\Delta+o(1)$ as, for every $\ell\in[k-1]$, a.a.s. there are no $\ell$-sets with more than $\Gamma_{\ell}$ common neighbors in $G(n,p)$. Therefore, $\delta_n=o(1)$. Thus, ${\sf P}(X_n=m-1)$ is equal to $h_{m-1}+o(1)$ where $h_{m-1}$ is the probability that the maximum $i$ such that there exist $i$ disjoint $k$-sets $U_1,\ldots,U_{i}\in{[n]\choose k}$ having $|N_n(U_j)|>b$ for all $j\in[i]$ equals $m-1$. But this probability is much easier to estimate. Indeed, set
$$
G(m-1)=\frac{1}{(m-1)!}\prod_{j=1}^{m-1}{n-k(j-1)\choose k},\quad
n_{m-1}=n-k(m-1),\quad
b_{m-1}=b-k(m-1).
$$
Fix disjoint $k$-sets $U_1,\ldots,U_{m-1}$. For a $k$-set $U$, denote by $\hat N_n(U)$ the number of common neighbors of $U$ in $\mathcal{N}_{m-1} = [n]\setminus[U_1\sqcup\ldots\sqcup U_{m-1}]$.
Then
$$
 G(m-1)\left({\sf P}(|\hat N_n(U_1)|>b)\right)^{m-1}{\sf P}\left(\Delta^1_{k;n_{m-1}}\leq b\right)\leq h_{m-1}\leq
$$
$$
 G(m-1)\left({\sf P}(|\hat N_n(U_1)|>b_{m-1}+k)\right)^{m-1}{\sf P}\left(\Delta^1_{k;n_{m-1}}\leq b_{m-1}\right).
$$
These bounds, in particular, follow from the fact that all $U_i$ are disjoint and the numbers of their common neighbors in $\mathcal{N}_{m-1}$ are not greater than $|N_n(U_i)|$ but not less than $|N_n(U_i)|-k(m-1)+k$. As $b_{m-1}=b+O(1)$ and $n_{m-1}=n+O(1)$, both lower and upper bounds are equal to $\frac{1}{(m-1)!}\left(e^{-y}\right)^{m-1}e^{-e^{-y}}(1+o(1))$. Therefore, reminding that $y$ and $m$ are constants, we finally derive
$$
 {\sf P}(\Delta^m_{k;n}\leq b)=\sum_{i=0}^{m-1}{\sf P}(X_n=i)\sim\sum_{i=0}^{m-1}h_i\sim e^{-e^{-y}}\sum_{i=0}^{m-1}\frac{1}{i!}\left(e^{-y}\right)^{i}.
$$

\section{Discussions and further questions}
\label{Further}

In Section~\ref{Intro}, we have mentioned that the dependencies between degrees of $G(n,p)$ are weak enough, and so, the result of Bollob\'{a}s does not essentially differ from the result of Nadarajah and Mitov (for independent binomial random variables). This motivates the following question. How strong can be dependencies between the binomial random variables until~(\ref{approaches_Gumbel}) fails?

Let us formalize this question in the following way. Let $\xi_1,\xi_2,\ldots$ be independent Bernoulli random variables with parameter $p$. For every $n\in\mathbb{N}$, consider an $M=M(n)$-element set $\Sigma_n$ of $n$-vectors $(\xi_{i_1},\ldots,\xi_{i_n})$ having components in the given sequence, $M\gg(\ln n)^3$. Let $D_n=\max_{(\xi_{i_1},\ldots,\xi_{i_n})\in\Sigma_n}(\xi_{i_1}+\ldots+\xi_{i_n})$.  Assume that $m=m(n)$ is such that any two vectors from $\Sigma_n$ have at most $m$ common components. For $m=0$, we have~(\ref{approaches_Gumbel}), and the result of Bollob\'{a}s relates to $m=1$. Can we guarantee the same for larger $m$?

Below, we state a generalization of our Janson-type inequality given in Section~\ref{Proof_Janson}, that has an analogous proof and immediately implies the following answer on the above question. If $m=o(\sqrt{n/\ln n})$, then (\ref{approaches_Gumbel}) is true.

\begin{lemma}
For every $n\in\mathbb{N}$, consider a sequence of independent random variables $\xi^n=(\xi^n_1,\xi^n_2,\ldots)$, an $M(n)$-element set $\{\eta^n_1,\ldots,\eta^n_M\}$ of vectors having components in $\xi^n$ and Borel sets $A^n_1,\ldots,A^n_M$, where $A^n_i\subset\mathbb{R}^{k_i}$, $k_i$ is the dimension of $\eta^n_i$. Let, for every $i\in\{1,\ldots,M\}$, $J(i)\subset\{1,\ldots,M\}$ be such that, for every $J\subseteq J(i)$,
$$
 {\sf P}(\wedge_{j\in J}\{\eta^n_j\notin A_j^n\}|\eta^n_i)={\sf P}(\wedge_{j\in J}\{\eta^n_j\notin A_j^n\})+o(1)\quad\text{a.s.}
$$
uniformly over all $i$ and $J$. Let
$$
X_n=\sum_{i=1}^M I(\eta^n_i\in A^n_i),\quad \Delta_n=\sum_{i=1}^n\sum_{j\in[M]\setminus J(i)}{\sf P}(\eta^n_i\in A^n_i,\eta^n_j\in A^n_j).
$$
If ${\sf P}(X_n=0)\geq e^{-{\sf E}X_n+o(1)}$, then
$$
{\sf P}(X_n=0)\leq\mathrm{exp}\left[-(1+o(1)){\sf E}X_n-(1+o(1))e^{{\sf E}X_n}\Delta_n\right].
$$
\end{lemma}




It is also of interest to consider the case $k=k(n)$ as well as $1-p$ approaching $0$ faster than $\sqrt{\frac{\ln\ln n}{\ln n}}$ and to prove (or disprove) an analogue of Theorem~\ref{main}.\\

Finally, in Section~\ref{Intro}, we mentioned that Theorem~\ref{main} can be formulated in terms of the maximum number of $(G,H)$-extensions where $|V(H)|-|V(G)|=1$ and $|E(H)|-|E(G)|=|V(G)|$. We expect that our techniques may give analogous results for a wider class of grounded strictly balanced pairs $(G,H)$.

\section{Acknowledgement}
The paper is partially supported by the Russian Foundation for Basic Research (grant 20-31-70025). The work of I. V. Rodionov in Sections 2.1, 2.3 was performed at the Institute for Information Transmission Problems (Kharkevich Institute) of the Russian Academy of Sciences with the support of the Russian Science Foundation (grant 21-71-00035).


\begin{thebibliography}{99}

\bibitem{Alon} N.~Alon, J.H.~Spencer, \emph{The Probabilistic Method}, Third Edition, John Wiley $\&$ Sons, 2008.

\bibitem{Anderson} C.W. Anderson, S.G. Coles, J. H\H{u}sler, {\it Maxima of Poisson-like variables and related triangular arrays}, Ann. Appl. Probab., {\bf7}:4 (1997), 953--971.

\bibitem{Beirlant} J. Beirlant, Y. Goegebeur, J. Segers, J.L. Teugels, {\it Statistics of Extremes: Theory and Applications}, John Wiley $\&$ Sons, 2004.

\bibitem{Bollobas} B. Bollob\'{a}s, {\it Random Graphs}, 2nd Edition, Cambridge University Press, 2001.

\bibitem{Bol_Degree} B. Bollob\'{a}s, {\it The distribution of the maximum degree of a random graph}, Discrete Mathematics, {\bf 32} (1980), 201--203.

\bibitem{Bol_Degree2} B. Bollob\'{a}s, {\it Degree sequences of random graphs}, Discrete Mathematics, {\bf 33} (1981), 1--19.

\bibitem{BopponaSpencer} R. Boppona, J. Spencer, {\it A useful elementary correlation inequality}, J. Combin. Theory Ser. A, {\bf 50} (1989), 305--307.

\bibitem{David} H.A. David, {\it Order Statistics}, New York: John Wiley \& Sons, 1970.


\bibitem{Sielenou} P. Sielenou Dkengne, N. Eckert, P. Naveau, {\it A limiting distribution for maxima of discrete stationary triangular arrays with an application to risk due to avalanches}, Extremes, {\bf 19}:1 (2016), 25--40.


\bibitem{Gnedenko} B. Gnedenko, {\it Sur la Distribution Limite du Terme Maximum d'une S\'{e}rie Al\'{e}atoire}, Annals of Mathematics, {\bf 44}  (1943), 423--453.

\bibitem{Haan} L. Haan, A. Ferreira, {\it Extreme Value Theory. An Introduction}, Springer, 2006.


\bibitem{Ivchenko} G.I. Ivchenko, {\it On the asymptotic behavior of degrees of vertices in a random graph}, Theory Probab. Appl., {\bf 18}:1 (1973), 188--195.

\bibitem{Janson} S. Janson, T. \L uczak, A. Rucinski, {\it Random Graphs}, New York, Wiley, 2000.

\bibitem{Kleitman} D.J. Kleitman, {\it Families of non-disjoint subsets}, J. Combinatorial Theory, {\bf 1} (1966), 153--155.

\bibitem{Leadbetter} M.R. Leadbetter, G. Lindgren, H. Rootz\'{e}n, {\it Extremes and Related Properties of Random Sequences and Processes}, New York: Springer Verlag, 1983.

\bibitem{Independent} S. Nadarajah, K. Mitov, {\it Asymptotics of Maxima of Discrete Random Variables}, Extremes, {\bf 5} (2002), 287--294.


\bibitem{Riordan} O. Riordan, L. Warnke, {\it The Janson inequalities for general up-sets}, Random Structures $\&$ Algorithms, {\bf 46}:2 (2015), 391--395.

\bibitem{Warnke_extensions} M. \v{S}ileikis, L. Warnke, {\it Counting extensions revisited}, 2019, https://arxiv.org/pdf/1911.03012.pdf.







\bibitem{SpencerSuen} J. Spencer, {\it A useful elementary correlation inequality, II}, J. Combin. Theory Ser. A, {\bf 84} (1998), 95--98.

\bibitem{Spencer_ext_counting} J.H. Spencer, {\it Counting extensions}, J. of Comb. Th. Ser A, {\bf 55} (1990), 247--255.






\end{thebibliography}
\end{document}